\documentclass[a4paper,11pt]{amsproc}
\usepackage{calrsfs, amsfonts, amssymb, mathrsfs, eucal}
\usepackage{braket}
\usepackage[all]{xy}

\def\Z{\mathbb{Z}}

\def\P{\mathbb{P}}

\def\O{\mathscr{O}}
\def\m{\mathfrak{m}}
\def\til#1{\widetilde{#1}}
\def\ovl#1{\overline{#1}}
\def\cl{\mathop{\mathrm{cl}}\nolimits}
\def\rig{\mathop{\mathrm{rig}}\nolimits}

\def\int{\mathop{\mathrm{int}}\nolimits}

\def\Spec{\mathop{\mathrm{Spec}}\nolimits}

\def\Spf{\mathop{\mathrm{Spf}}\nolimits}
\def\Sph{\mathop{\mathrm{Sph}}\nolimits}

\def\Proj{\mathop{\mathrm{Proj}}\nolimits}
\def\Frac{\mathop{\mathrm{Frac}}\nolimits}

\def\Hom{\mathop{\mathrm{Hom}}\nolimits}

\def\lHom{\mathop{\mathscr{H}\!\mathit{om}}\nolimits}

\def\H{\mathrm{H}}

\def\an{\mathrm{an}}

\def\Itor{\textrm{$I$-}\mathrm{tor}}

\def\opp{\mathrm{opp}}
\def\dl{\langle\!\langle}
\def\dr{\rangle\!\rangle}

\def\Mod{\mathbf{Mod}}

\def\Sch{\mathbf{Sch}}

\def\Hs{\mathbf{Hs}}

\def\CHs{\mathbf{CHs}}

\def\CRh{\mathbf{CRh}}

\def\het#1{{#1}^h}

\def\ZR#1{\langle #1\rangle}

\def\Emb{\mathbf{Emb}}

\newcommand*{\longhookrightarrow}{\ensuremath{\lhook\joinrel\relbar\joinrel\rightarrow}}
\title[On Henselian Rigid Geometry]{On Henselian Rigid Geometry}
\author[F.\ Kato]{Fumiharu Kato}
\thanks{}
\begin{document}
\newtheorem*{thm*}{Theorem}
\newtheorem{thm}{Theorem}
\newtheorem{prop}[thm]{Proposition}
\newtheorem{lem}[thm]{Lemma}
\newtheorem*{cor*}{Corollary}
\newtheorem{cor}[thm]{Corollary}
\newtheorem{rem}[thm]{Remark}
\newtheorem{dfn}[thm]{Definition}
\newtheorem{ntn}[thm]{Notation}
\newtheorem{exa}[thm]{Example}
\newtheorem{exas}[thm]{Examples}
\newtheorem{note}[thm]{Note}
\newtheorem{probs}[thm]{Problems}
\newtheorem{assum}[thm]{Assumption}
\newtheorem{sit}[thm]{Situation}
\newtheorem{cla}[thm]{Claim}
\newtheorem{para}[thm]{}
\newtheorem{propsec}{Proposition}[section]
\maketitle
\begin{abstract}
We overview some of the foundations of the so-called henselian rigid geometry, and show that henselian rigid geometry has many aspects, useful in applications, that one cannot expect in the usual rigid geometry. This is done by announcing a few characteristic results, one of which is an analogue of Zariski Main Theorem.
\end{abstract}

\setcounter{tocdepth}{2}
\tableofcontents

\section{Introduction}\label{sec-introduction}
By the celebrated theorem by Raynaud, rigid geometry was suggested to have the characteristic architecture, built upon what we call a {\em model geometry}, a suitable framework of geometry equipped with a class of morphisms, called {\em admissible blow-ups}, from which the rigid geometry is induced by inverting all admissible blow-ups.
Here, one of the novelties is the viewpoint that rigid geometry can be regarded as the ``birational geometry'' of the model geometry, which allows one to envisage the so-called {\em birational approach} to the rigid geometry, as pursued in our book \cite{FK2}.
Another novel point that arises from this doctrine is that rigid geometry, as a whole, depends on the choice of the model geometry, for which, while the usual choice being the geometry of formal schemes, one can take fairly freely from diversity of geometries that appear in arithmetic or algebraic geometry.
Other than formal schemes, a promising candidate for a model geometry is the geometry of henselian schemes, from which one obtains the so-called {\em henselian rigid geometry}.

It is expected that, at least in the foundational level, henselian rigid geometry in many respects should go parallel to the usual rigid geometry without major modification; also expected is that henselian rigid geometry, in many situations, should be technically easier to develop than usual rigid geometry, due in large part to the fact that henselian rings and henselization of rings are technically easier to handle than complete rings and completion. 
However, although the first appearance of henselian rigid geometry in literature seems in \cite{Fujiw1}, henselian rigid geometry itself, since then and therein, seems not to have been studied well as an independent and autonomous discipline, but rather treated only as a passable substitute for the usual rigid geometry.

There should be several reasons for the sluggish pace of acceptance and development of henselian rigid geometry.
Perhaps the most important among them is that there have been no clues how henselian rigid geometry is different from usual rigid geometry, and, consequently, how it is especially useful in applications. 
Another reason, which we should point out here, is that the geometry of henselian schemes, which henselian rigid geometry adopt as its model geometry, has not yet been fully developed, although there are several fairly systematic accounts on it, e.g., \cite{Cox1}\cite{KPR}, especially those from italian school \cite{Greco1}\cite{GS}\cite{M}\cite{S}.
What is still missing is the henselian GAGA theory, or so to speak, GHGA ($=$ g\'eom\'etrie henselienne et g\'eom\'etrie alg\'ebrique), which should give us the comparison of the theories of coherent sheaves, comparison of cohomologies and existence (algebrizability), between henselian scheme geometry and the geometry of schemes, giving thereby the basis for GAGA between henselian rigid geometry and algebraic geometry.
In fact, it seems that GHGA is a very difficult theorem, and we still do not know if the full version of it is true or not, whereas we have some partial results (as indicated in \S\ref{sub-GHGA}).

This is the proceedings report of what the author has delivered in his talk at Algebraic Geometry in East Asia 2016.
In the talk, the author has reported his recent joint-work with Shuji Saito (Tokyo Institute of Technology), the main result of which is an analogue of Zariski Main Theorem in henselian rigid geometry, which, we think, together with a few of its corollaries, very clearly illustrates how henselian rigid geometry is different from usualy rigid geometry, and moreover, convinces us that henselian geometry is indeed useful, in a way that classical theory could not afford, for application to algebraic geometry.

Our ZMT theorem gives, similarly to the classical one in algebraic geometry, an assertion of the following sort: a quasi-finite mapping factorizes by an open immersion followed by a finite map.
Here, in our case, while the domain of the quasi-finite map can be an arbitrary henselian affinoid space, the target has to be the analytification of an algebraic variety.
One can then show that the assertion of the above-mentioned type is true, where the finite map in the factorization can be chosen to be algebraic.

Now let us describe our results.
Let $K$ be a non-trivially valued henselian valuation field, i.e., the fractional field of an $a$-adically henselian valuation ring $V$, where $a$ is a non-zero element in the maximal ideal $\m_V$ of $V$.
Let $X$ be a separated of finite type scheme over $K$, and consider its analytification $X^{\an}$, which is a separated locally of finite type henselian rigid space over $K$.

\begin{thm*}[{\sc $=$ Theorem \ref{thm-ZMT}}]\label{thm-I}
Let $\varphi\colon\mathscr{U}\hookrightarrow X^{\an}$ be a quasi-finite $K$-map from a henselian affinoid space $\mathscr{U}$ of finite type over $K$.
Then there exists a finite morphism $g\colon W\rightarrow X$ with the diagram
$$
\xymatrix{\mathscr{U}\ar[r]_{j}\ar@/^1pc/[rr]^{\varphi}&W^{\an}\ar[r]_{g^{\an}}&X^{\an},}
$$
where $j$ is an open immersion.
\end{thm*}

From the theorem and the techniques of the proof, one has several useful consequences; the following is one of them:

\begin{cor*}[{\sc $=$ Corollary \ref{cor-ZMT}}]\label{cor-I}
Let $X$ be a separated finite type scheme over $K$, $\mathscr{U}$ a henselian rigid space of finite type over $K$, and $\varphi\colon \mathscr{U}\hookrightarrow X^{\an}$ an immersion.
Then there exists a {\em closed subscheme} $W\subseteq X$ that is smallest among those containing the image of $\mathscr{U}$ as an open subspace.
\end{cor*}

The point lies in that, in the corollary, the closed subscheme $W\subseteq X$ gives the {\em scheme-theoretic closure} of the affinoid subspace $\mathscr{U}$ in $X$.
Notice that these kinds of results, ZMT and scheme-theoretic closure, are far from being true in the usual rigid geometry, since, for example, analytic subspaces can be highly transcendental.

It may be that, in view of the facts that ``henselian'' is like ``algebraic'', and that henselization is like algebraic closure, these results are not extremely surprising.
But their significance lies rather in the fact that several algebraic constructions in algebraic geometry can be done analytic-locally first, and then, extended to the ``closure'', as indicated the above statements.
In fact, these results were originally intended to apply to the theory of algebraic cycles, in which Shuji Saito, by the joint-work with M.\ Kerz, tries to give a construction of what they call {\em analytic Chow groups}.
In this context, in particular, the following theorem, which is also among the corollaries of our technique, should be important:
Let $X$ be a proper scheme over $K$, and $\mathscr{U}$ (resp.\ $U$) a finite type henselian rigid space (resp.\ finite type scheme) over $K$.
As usual, a flat family of closed subspaces in $X$ over $\mathscr{U}$ (resp.\ $U$) is a closed subspace $\mathscr{Y}\subseteq X^{\an}\times_K\mathscr{U}$ (resp.\ $Y\subseteq X\times_KU$) that is flat over $\mathscr{U}$ (resp.\ $U$).

\begin{thm*}[{\sc $=$ Theorem \ref{thm-FCS}}]\label{thm-II}
For any flat family $\mathscr{Y}\subseteq X^{\an}\times_K\mathscr{U}$ over finite type henselian affinoid $\mathscr{U}$, there exists an affine finite type scheme $U$ such that
\begin{itemize}
\item[{\rm (a)}] $U^{\an}$ contains $\mathscr{U}$ as an affinoid subdomain$;$
\item[{\rm (b)}] $\mathscr{Y}$ extends to a flat family $Y\subseteq X\times_KU$ over $U$.
\end{itemize}
\end{thm*}

The proof of the above theorems, as well as the statements themselves, partly depend on henselian rigid GAGA theorems.
As mentioned above, these theories are still missing in literature, and so one has to first work out these foundations, before coming up to the proof of the above theorems.
As this article is intended to give a first announcement of our results, we are only to give a brief survey of our theory.
The more precise and systematic accounts will come elsewhere in the following form: The foundational part, consisting of generalities of the geometry of henselian schemes, affineness criterion, cohomology calculus, and GHGA theorems, will be written by the author, and the henselian rigid ZMT, as well as the related results, will be written by collaboration of Shuji Saito and the author.  

Under these circumstances, the most reasonable role of this article is to give a survey of these forthcoming papers.
So, in the next section, we will give a survey of the first paper, the one on the foundations of henselian scheme geometry, including some directions to the part of GHGA theorems, which we can so far prove.
Then, in section \S\ref{sec-henselianrigidspaces}, we describe roughly the contents of what we will present in the second paper, including the sketches of the proofs of the above results.

As already mentioned, this article is the proceedings report of my talk in the conference ``Algebraic Geometry in East Asia 2016'' in January 2016. 
The author thanks the organizers of the conference for the invitation.
Thanks are also due to the referee for the valuable comments and suggestions, which fill in gaps and fix errors in the first draft of this paper.

\section{Henselian schemes}\label{sec-henselianschemes}
\subsection{Henselian finite type algebras}\label{sub-topfintype}
We refer to \cite{Cox1}\cite{Greco1}\cite{KPR} for the first generalities on henselian rings and henselian schemes.
For a ring $A$ with an ideal $I\subseteq A$, the {\em $I$-adic henselization} is denoted by $\het{A}$.
The natural map $A\rightarrow\het{A}$ is always flat, and is faithfully flat if $I$-adically Zariskian, i.e., $1+I\subseteq A^{\times}$ (cf.\ \cite[{\bf 0}.7.3.8]{FK2}).

Let $A$ be an $I$-adically henselian ring by an ideal $I\subseteq A$.
We denote by
$$
A\{X_1,\ldots,X_n\}
$$
the $I$-adic henselization of the polynomial ring $A[X_1,\ldots,X_n]$.
We say that an $A$-algebra $B$ is {\em of $($henselian$)$ finite type} (resp.\ {\em of $($henselian$)$ finite presentation}) if $B$ is isomorphic to the $A$-algebra of the form $A\{X_1,\ldots,X_n\}/\mathfrak{a}$ by an ideal (resp.\ finitely generated ideal) $\mathfrak{a}$ of $A\{X_1,\ldots,X_n\}$.
Notice that these $A$-algebras are again $I$-adically henselian.
One can prove that the $I$-adic henselization of a finite type $A$-algebra is of henselian finite type.
Moreover, an $A$-algebra is of henselian finite presentation if and only if it is the $I$-adic henselization of a finitely presented $A$-algebra.

Recall that an $I$-adically topologized ring $A$ of finite ideal type is said to be {\em $I$-adically adhesive} if the following conditions are satisfied (\cite[Chap.\ {\bf 0}, \S8.5.(a)]{FK2}):
\begin{itemize}
\item[{\rm (a)}] $\Spec A\setminus V(I)$ is a Noetherian scheme$;$
\item[{\rm (b)}] for any finitely generated $A$-module $M$, the $I$-torsion part $M_{\Itor}$ is finitely generated.
\end{itemize}
If, moreover, any finite type $A$-algebra is $I$-adically adhesive, $A$ is said to be {\em $I$-adically universally adhesive}.

It is known (\cite[Chap.\ {\bf 0}, \S8.5.(c)]{FK2}) that $I$-adically adhesive ring $A$ has the {\it adicness-preservation} property (denoted by {\bf (AP)} in [loc.\ cit.]), i.e., for any finitely generated $A$-module $M$ and any $A$-submodule $N\subseteq M$, the subspace topology on $N$ induced from the $I$-adic topology on $M$ is the $I$-adic topology on $N$.
From this, it follows that $I$-adically adhesive rings enjoy several pleasant properties; for example,
\begin{itemize}
\item the functor $M\mapsto\widehat{M}$ by $I$-adic completion on the full subcategory of $\Mod_A$ consisting of finitely generated $A$-modules is exact;
\item $M\otimes_A\widehat{A}\cong\widehat{M}$ for any finitely generated $A$-module $M$;
\item the $I$-adic completion map $A\rightarrow\widehat{A}$ is flat.
\end{itemize}
If, moreover, $A$ is $I$-adically henselian, then (\cite[Chap.\ {\bf 0}, (7.4.16)]{FK2}):
\begin{itemize}
\item any finitely generated $A$-module is $I$-adically separated; in particular, $A$ itself is automatically $I$-adically separated;
\item any $A$-submodule $N$ of a finitely generated $A$-module $M$ is closed in $M$ with respect to the $I$-adic topology.
\end{itemize}
Here is one more good thing about adhesiveness:
\begin{itemize}
\item if $A$ is $I$-adically universally adhesive, and $I$-torsion free, then any finitely presented $A$-algebra $B$ is a coherent ring, and hence, any finitely presented $B$-module is a coherent module.
\end{itemize}
(Recall that a module $M$ over a ring $A$ is said to be {\it coherent} if it is finitely generated, and every finitely generated $A$-submodule is finitely presented, and that $A$ is a {\it coherent ring} if it is coherent as a module over itself, i.e., every finitely generated ideal is finitely presented; cf.\ \cite[Chap.\ I, \S 3, Exercise 11]{Bourb1}.)

\begin{thm}\label{thm-henselianvaluationrings0}
Let $V$ be an $a$-adically henselian valuation ring, where $a\in\m_V\setminus\{0\}$.
Then any henselian finite type $V$-algebra is $a$-adically universally adhesive.
\end{thm}

To show that $V$ is $a$-adically universally adhesive, we only need to show that $A=V[X_1,\ldots,X_n]$ is $a$-adically adhesive (cf.\ \cite[Chap.\ {\bf 0}, (8.5.7) (2)]{FK2}).
Since $V$ is $a$-adically adhesive (\cite[Chap.\ {\bf 0}, (8.5.15)]{FK2}), $V\rightarrow\widehat{V}$ is flat (\cite[Chap.\ {\bf 0}, (8.2.18)]{FK2}), where $\widehat{V}$ is the $a$-adic completion of $V$.
Hence $V[X_1,\ldots,X_n]\rightarrow\widehat{V}[X_1,\ldots,X_n]$ is flat.
By \cite[Chap.\ {\bf 0}, (9.2.7)]{FK2}, the map $\widehat{V}[X_1,\ldots,X_n]\rightarrow\widehat{V}\dl X_1,\ldots,X_n\dr$ is flat.
Hence $A=V[X_1,\ldots,X_n]\rightarrow\widehat{A}=\widehat{V}\dl X_1,\ldots,X_n\dr$ is flat.
Then by \cite[Chap.\ {\bf 0}, (9.2.7), (8.5.11)]{FK2}, we deduce that $A$ is $a$-adically adhesive.

Next we show that the $a$-adic henselization $\het{A}=V\{X_1,\ldots,X_n\}$ is $a$-adically adhesive.
Take a filtered inductive system $\{A_{\lambda}\}_{\lambda\in\Lambda}$ consisting of finite type $V$-algebras with \'{e}tale maps such that $\het{A}=\varinjlim_{\lambda\in\Lambda}A_{\lambda}$ and $\het{(A_{\lambda})}\cong\het{A}$.
We already know that each $A_{\lambda}$ is $a$-adically adhesive, and hence $A_{\lambda}\rightarrow\widehat{A}_{\lambda}=\widehat{A}$ is flat.
For any finitely generated ideal $J\subseteq\het{A}$, take sufficiently large $\lambda$ and a finitely generated ideal $J_{\lambda}\subseteq A_{\lambda}$ such that $J=J_{\lambda}\het{A}=J_{\lambda}\otimes_{A_{\lambda}}\het{A}$.
Since $\widehat{A}$ is flat over $A_{\lambda}$, we have $J\otimes_{\het{A}}\widehat{A}=J_{\lambda}\otimes_{A_{\lambda}}\widehat{A}=J_{\lambda}\widehat{A}=J\widehat{A}$.
Thus we deduce that $\het{A}\rightarrow\widehat{A}$ is flat, in particular, faithfully flat.
Now, since $\widehat{A}[\frac{1}{a}]$ is known to be Noetherian (\cite[Chap.\ {\bf 0}, (9.2.7)]{FK2}), we deduce that $\het{A}[\frac{1}{a}]$ is Noetherian.
By \cite[Chap.\ {\bf 0}, (9.2.7), (8.5.11)]{FK2}, we conclude that $\het{A}$ is $a$-adically adhesive, as desired.

\subsection{Henselian schemes}\label{sub-henselianschemes}
Let $A$ be an $I$-adically henselian ring by an ideal $I\subseteq A$.
The {\em henselian spectrum} of $A$, denoted by $X=\Sph A$, is the topologically locally ringed space defined as follows (cf.\ \cite[1.8, 1.9]{GS}):
\begin{itemize}
\item it is, as a set, the subset $V(I)$ of $\Spec A$, or equivalently, the set of all open prime ideals of $A$;
\item the topology is the subspace topology induced from the Zariski topology of $\Spec A$;
\item for any $f\in A$, set $\mathcal{D}(f)=D(f)\cap X$ and $\O_X(\mathcal{D}(f))=\het{A}_f$; then this defines a sheaf $\O_X$ of rings on $X$.
\end{itemize}
Notice that, for any $x=\mathfrak{p}\in X$, the stalk $\O_{X,x}$ is given by $\het{A}_{\mathfrak{p}}$, the $I$-adic henselization of $A_{\mathfrak{p}}$.
The open subset $\mathcal{D}(f)=D(F)\cap X$, considered as an open subspace of $X$, is isomorphic to the henselian spectrum $\Sph\het{A}_f$.

\begin{dfn}\label{dfn-henselianschemes}{\rm 
(1) An {\em affine henselian scheme} is a topologically locally ringed space isomorphic to $(X=\Sph A,\O_X)$ for an $I$-adically henselian ring $A$.

(2) A topologically locally ringed space $(X,\O_X)$ is called a {\em henselian scheme} if it has an open covering $X=\bigcup_{\alpha\in L}U_{\alpha}$ by affine henselian schemes.}
\end{dfn}

A morphism between henselian schemes is a morphism of topologically locally ringed spaces.
Any morphism $f\colon\Sph B\rightarrow\Sph A$, where $A$ (resp.\ $B$) is considered with the $I$-adic (resp.\ $J$-adic) topology, comes uniquely from a continuous homomorphism $\varphi\colon A\rightarrow B$, where the continuity is equivalent to that there exists $n\geq 0$ such that $I^nB\subseteq J$; viz., the functor $A\mapsto\Sph A$ gives the categorical equivalence between the opposite category of the category of henselian rings with continuous homomorphisms to the category of affine henselian schemes.
If, in the above situation, $IB$ gives an ideal of definition of the topological ring $B$, i.e., $IB$-adic topology coincides with the $J$-adic topology, then we say that the morphisms $f$ and $\varphi$ are {\em adic}.
Adic morphisms between general henselian schemes are defined in the obvious manner.

\subsection{Henselization of schemes}\label{sub-henselizationscheme}
Let $X$ be a scheme, and $Y\subseteq X$ a closed subscheme.
The {\em henselization of $X$ along $Y$}, denoted by $\het{X}|_Y$, or by $\het{X}$, is the henselian scheme defined as follows (cf.\ \cite[I, \S1]{Cox1}):
\begin{itemize}
\item the underlying topological space is given by the underlying topological space of $Y$;
\item for any affine open subset $U\cong\Spec A$ of $X$, $\O_{\het{X}}(U\cap Y)$ is given by the $I$-adic henselization $\het{A}$ of $A$, where $I\subseteq A$ is the defining ideal of $Y\cap U$ in $U$.
\end{itemize}

Notice that there is a natural morphism $j\colon\het{X}|_Y\rightarrow X$ of locally ringed spaces, which is flat in the sense that, for any $x\in\het{X}$, the morphism $\O_{X,j(x)}\rightarrow\O_{\het{X},x}$ is flat.

\subsection{Quasi-coherent and coherent sheaves}\label{sub-quasicoherentsheaves}
The theory of quasi-coherent sheaves on henselian schemes has already been fully developed in \cite{GS}.
Let $X=\Sph A$ be an affine henselian scheme. 
$X$ can be viewed as the henselization of $Y=\Spec A$ along the closed subscheme $V(I)$ corresponding to an ideal of definition $I\subseteq A$.
Let $j\colon X\rightarrow Y$ be the canonical morphism of locally ringed spaces.
For any $A$-module $M$, the pull-back $j^{\ast}\til{M}$, which we again denote simply by $\til{M}$, of the quasi-coherent sheaf $\til{M}$ on $Y$ is a quasi-coherent sheaf on $X$.
As in the scheme case, the functor by $M\mapsto\til{M}$ is an exact equivalence between the category of $A$-modules and the category of quasi-coherent sheaves on $X$, with the quasi-inverse given by $\mathscr{F}\mapsto\Gamma_X(\mathscr{F})$.
Moreover, we have the following analogue of ``Theorem B'': 
\begin{thm}[{\rm \cite[1.12]{GS}}]
For any quasi-coherent sheaf $\mathscr{F}$ on $X=\Sph A$, we have $\H^q(X,\mathscr{F})=0$ for $q\geq 1$.
\end{thm}

Let $V$ be an $a$-adically henselian valuation ring, where $a$ is a non-zero element in the maximal ideal $\m_V$ of $V$.
It follows from Theorem \ref{thm-henselianvaluationrings0} that, if $X$ is a henselian scheme flat of finite type over $\Sph V$, then the structure sheaf $\O_X$ is coherent.
In particular, an $\O_X$-module $\mathscr{F}$ is coherent if and only if it is of finite presentation.
In case $X$ is affine $X=\Sph A$, where $A$ is an $a$-torsion free henselian of finite type $V$-algebra, then $M\mapsto\til{M}$ gives the exact equivalence between the category of coherent (equivalently, finitely presented) $A$-modules to the category of coherent sheaves on $X$.

\subsection{GHGA}\label{sub-GHGA}
Let $Y$ be a proper and flat scheme over $\Spec V$, where $V$ is an $a$-adically henselian valuation ring, and $X=\het{Y}$ the $a$-adic henselization, which is a proper and flat henselian scheme over $\Sph V$.
We have moreover the $a$-adic completion $\widehat{X}$, which is a proper and flat formal scheme over $\Spf\widehat{V}$.
Notice that the completion map $V\rightarrow\widehat{V}$ is fathfully flat.
Consider the maps
$$
\widehat{X}\stackrel{\phi}{\longrightarrow}X\stackrel{\varphi}{\longrightarrow}Y
$$
of locally ringed spaces, which are flat.
We have the chain of exact functors
$$
\Mod_Y\stackrel{\varphi^{\ast}}{\longrightarrow}\Mod_X\stackrel{\phi^{\ast}}{\longrightarrow}\Mod_{\widehat{X}},
$$
which map coherent sheaves to coherent sheaves.

The GHGA ($=$ henselian geometry and algebraic geometry) statements that we can state here are as follows.
\begin{thm}[GHGA comparison for $\H^0$]\label{thm-H0}
For any coherent sheaf $\mathscr{F}$ on $Y$, the canonical morphism
$$
\H^0(Y,\mathscr{F})\longrightarrow\H^0(X,\varphi^{\ast}\mathscr{F})
$$
is an isomorphism.
\end{thm}

\begin{thm}[GHGA existence for subquotients]\label{thm-subquotients}
For any coherent sheaf $\mathscr{F}$ on $Y$, any coherent subsheaf $\mathscr{G}$ of $\varphi^{\ast}\mathscr{F}$ is algebrizable, i.e., $\mathscr{G}=\varphi^{\ast}\mathscr{H}$ for a coherent subsheaf $\mathscr{H}$ of $\mathscr{F}$.
\end{thm}

Let us show Theorem \ref{thm-H0}. We first claim the following.

\medskip
{\sc Claim 1.} {\it Let $A$ be a $V$-algebra of finite type, and $Y$ a proper scheme over $A$.
Set $X=\het{Y}$, which is a proper henselian scheme over $\Sph\het{A}$.
Let $\varphi\colon X\rightarrow Y$ be the canonical map.
Then, for any coherent sheaf $\mathscr{F}$ on $Y$, $\H^0(X,\varphi^{\ast}\mathscr{F})$ is a coherent $A$-module.}

\medskip
We first show the claim in case $Y$ is projective over $A$.
Take a closed immersion $\iota\colon Y\hookrightarrow\P^n_A=\Proj S$, where $S=A[x_0,\ldots,x_n]$.
With suitable choice of $n$ and the homogeneous coordinates, we may assume that $\iota_{\ast}\mathscr{F}$ is $x_i$-torsion free for $i=0,\ldots,n$.
Then, considering $\iota_{\ast}\mathscr{F}$ in place of $\mathscr{F}$, we may assume $Y=\P^n_A$ and that $\mathscr{F}$ is $x_i$-torsion free for $i=0,\ldots,n$.

Set $M=\bigoplus_{k\in\Z}\Gamma(Y,\mathscr{F}(k))$, which is $x_i$-torsion free for $i=0,\ldots,n$.
Let $\{U_i=\Spec S_{(x_i)}:i=0,\ldots,n\}$ be the standard affine covering of $Y$.
Then $\{\het{U}_i=\Sph \het{(S_{(x_i)})}:i=0,\ldots,n\}$ gives an affine covering of $X$.
Since $\varphi^{\ast}\mathscr{F}$ is coherent on $X$, we have $\Gamma(\het{U}_i,\varphi^{\ast}\mathscr{F})=\Gamma(U_i,\mathscr{F})\otimes_{S_{(x_i)}}\het{(S_{(x_i)})}$, which coincides with the homogenous degree-zero part of $M\otimes_S\het{(S_{x_i})}$ for each $i=0,\ldots,n$; note that $\het{(S_{x_i})}=\het{A}\{x_0,\ldots,x_n,x^{-1}_i\}$ has, viewed as a subring of formal power series rings, the natural notion of degree for monomials, and the degree-zero part of $\het{(S_{x_i})}$ is a closed (hence $a$-adically henselian) subring, which thereby coincides with $\het{(S_{(x_i)})}$.
Hence elements of $\H^0(X,\varphi^{\ast}\mathscr{F})$ are written as $n+1$-tuples $(s_0,\ldots,s_n)$ with each homogenous $s_i$ of degree zero in $M\otimes_S\het{(S_{x_i})}$ such that $s_i$ and $s_j$ coincide in $M\otimes_S\het{(S_{x_ix_j})}$.
Since $M$ is $x_i$-torsion free, these modules are submodules in $M\otimes_S\het{A}\{x^{\pm 1}_0,\ldots,x^{\pm 1}_n\}$.
Then, similarly to the classical argument for calculating cohomologies over the projective space, one shows that $\H^0(X,\varphi^{\ast}\mathscr{F})\cong M_0\otimes_A\het{A}=\H^0(Y,\mathscr{F})\otimes_A\het{A}$, which is a coherent $V$-module.

To proceed, we need to show the following.

\medskip
{\sc Claim 2.} {\it Let $\pi\colon\til{Y}\rightarrow Y$ be a projective morphism between $V$-schemes of finite type, and set $X=\het{Y}$ and $\til{X}=\het{\til{Y}}$.
Consider the commutative diagram
$$
\xymatrix{\til{X}\ar[r]^{\til{\varphi}}\ar[d]_{\eta}&\til{Y}\ar[d]^{\pi}\\ X\ar[r]_{\varphi}&Y,}
$$
where the horizontal arrows are the canonical ones, and $\eta=\het{\pi}$.
Then, for any coherent sheaf $\mathscr{G}$ on $\til{Y}$, we have $\eta_{\ast}\til{\varphi}^{\ast}\mathscr{G}\cong\varphi^{\ast}\pi_{\ast}\mathscr{G}$.}

\medskip
To see this, it suffices to show that, for any affine open subset $U=\Spec A$ of $Y$, where $A$ is a $V$-algebra of finite type, $\H^0(\eta^{-1}\varphi^{-1}(U),\til{\varphi}^{\ast}\mathscr{G})\cong \H^0(\pi^{-1}(U),\mathscr{G})\otimes_A\het{A}$.
Hence we may assume $Y=\Spec A$, and what we need to show is that $\H^0(\til{X},\til{\varphi}^{\ast}\mathscr{G})\cong\H^0(\til{Y},\mathscr{G})\otimes_A\het{A}$, which we have already shown above.

Now, to show Claim 1 for general $Y$, by Carving lemma (\cite[Chap.\ {\bf I}, (8.3.2)]{FK2}), one can reduce to the following situation: There exists a closed immersion $\iota\colon Y_1\hookrightarrow Y$ and a projective morphism $\pi\colon\til{Y}\rightarrow Y$, which is isomorphic over $Y\setminus Y_1$, such that the claim is true for $Y_1$ and that $\til{Y}$ is projective over $V$.
In this situation, let $\mathscr{N}$ be the cokernel of $\mathscr{F}\hookrightarrow\pi_{\ast}\pi^{\ast}\mathscr{F}$.
Since the claim is true for the projective $\til{Y}$, we already know that $\H^0(X,\varphi^{\ast}\pi_{\ast}\pi^{\ast}\mathscr{F})$, which is isomorphic to $\H^0(\til{X},\til{\varphi}^{\ast}\pi^{\ast}\mathscr{F})$ due to Claim 2, is a coherent $V$-module.
Moreover, since $\mathscr{N}$ is coherent on $Y_1$, we also know that $\H^0(X,\varphi^{\ast}\mathscr{N})$ is coherent by our assumption.
Having the exact sequence
$$
0\longrightarrow\varphi^{\ast}\mathscr{F}\longrightarrow\varphi^{\ast}\pi_{\ast}\pi^{\ast}\mathscr{F}\longrightarrow\varphi^{\ast}\mathscr{N}\longrightarrow 0
$$
(due to the flatness of $\varphi$), we deduce the kernel $\H^0(X,\varphi^{\ast}\mathscr{F})$ of the map $\H^0(X,\varphi^{\ast}\pi_{\ast}\pi^{\ast}\mathscr{F})\rightarrow\H^0(X,\varphi^{\ast}\mathscr{N})$ is also coherent, which finishes the proof of Claim 1.

Now, to show Theorem {\rm \ref{thm-H0}}, take a finite affine covering $\mathscr{U}=\{U_{\alpha}=\Spec A_{\alpha}\}$ of $Y$, which induces affine coverings $\het{\mathscr{U}}=\{\het{U}_{\alpha}=\Spec\het{A}_{\alpha}\}$ and $\widehat{\mathscr{U}}=\{\widehat{U}_{\alpha}=\Spec\widehat{A}_{\alpha}\}$ of $X$ and $\widehat{X}$, respectively.
Set $M_{\alpha_0\cdots\alpha_p}=\Gamma(U_{\alpha_0\cdots\alpha_p},\mathscr{F})$ (where $U_{\alpha_0\cdots\alpha_p}=U_{\alpha_0}\cap\cdots\cap U_{\alpha_p}$).
Then the \v{C}ech modules $\mathscr{C}^p(\het{\mathscr{U}},\varphi^{\ast}\mathscr{F})$ is given by
$$
\mathscr{C}^p(\het{\mathscr{U}},\varphi^{\ast}\mathscr{F})=\prod_{\alpha_0,\ldots,\alpha_p}M_{\alpha_0\cdots\alpha_p}\otimes_{A_{\alpha_0\cdots\alpha_p}}\het{A}_{\alpha_0\cdots\alpha_p}, 
$$
and its $a$-adic completion coincides, up to canonical isomorphism, with $\mathscr{C}^p(\widehat{\mathscr{U}},\phi^{\ast}\varphi^{\ast}\mathscr{F})$.
Set $C^p=\mathscr{C}^p(\het{\mathscr{U}},\varphi^{\ast}\mathscr{F})$ and $B^p=\mathscr{C}^p(\het{\mathscr{U}},\O_X)$ for $p\geq 0$.
$K^0=\H^0(X,\varphi^{\ast}\mathscr{F})$ is the kernel of $C^0\rightarrow C^1$, and $\H^0(\widehat{X},\phi^{\ast}\varphi^{\ast}\mathscr{F})$ is the kernel of $\widehat{C}^0\rightarrow\widehat{C}^1$, where $\widehat{C}^p$ denotes the $a$-adic completion of $C^p$.
We have an injection $\varprojlim_{n\geq 0}K^0/K^0\cap a^nC^0\hookrightarrow\widehat{C}^0$.
Now, since $C^0$ is finitely generated over $B^0$, which is a henselian of finite type (hence $a$-adically adhesive) $V$-algebra, and since $K^0\subseteq C^0$ is a submodule over $B^0$, we deduce that the topology on $K^0$ by the filtration $\{K^0\cap a^nC^0\}$ is $a$-adic, and hence that $\varprojlim_{n\geq 0}K^0/K^0\cap a^nC^0\cong\widehat{K}^0$, the $a$-adic completion of $K^0$.
This shows that $\widehat{K}^0\rightarrow\widehat{C}^0$ is injective, and hence $\widehat{K}^0\rightarrow\H^0(\widehat{X},\phi^{\ast}\varphi^{\ast}\mathscr{F})$ is also injective; note that, as we already know by GFGA, $\H^0(\widehat{X},\phi^{\ast}\varphi^{\ast}\mathscr{F})$ is a finitely generated (hence $a$-adically complete) $\widehat{V}$-module.
Since we have $\H^0(X,\varphi^{\ast}\mathscr{F})\otimes_V\widehat{V}\cong\widehat{K}^0$ due to Claim, we deduce that the morphism $\H^0(X,\varphi^{\ast}\mathscr{F})\otimes_V\widehat{V}\rightarrow\H^0(\widehat{X},\phi^{\ast}\varphi^{\ast}\mathscr{F})$ is injective.

Now, by GFGA comparison, the composition
$$
\H^0(Y,\mathscr{F})\otimes_V\widehat{V}\longrightarrow\H^0(X,\varphi^{\ast}\mathscr{F})\otimes_V\widehat{V}\longhookrightarrow\H^0(\widehat{X},\phi^{\ast}\varphi^{\ast}\mathscr{F})
$$
is an isomorphism, and hence the morphisms therein are all isomorphisms.
Since $\widehat{V}$ is faithfully flat over $V$, we have $\H^0(Y,\mathscr{F})\cong\H^0(X,\varphi^{\ast}\mathscr{F})$, which finishes the proof of Theorem {\rm \ref{thm-H0}}.

To show Theorem \ref{thm-subquotients}, note that, by GFGA, we have $\mathscr{H}$ such that $\phi^{\ast}\mathscr{G}\cong\phi^{\ast}\varphi^{\ast}\mathscr{H}$.
Hence we need to show: For a coherent sheaf $\mathscr{F}$ on $X$ and coherent subsheaves $\mathscr{G}_1,\mathscr{G}_2\subseteq\mathscr{F}$, $\phi^{\ast}\mathscr{G}_1=\phi^{\ast}\mathscr{G}_2$ as a coherent subsheaf of $\phi^{\ast}\mathscr{F}$ implies $\mathscr{G}_1=\mathscr{G}_2$.
The condition reads $\mathscr{G}_1/\mathscr{G}_1\cap a^n\mathscr{F}=\mathscr{G}_2/\mathscr{G}_2\cap a^n\mathscr{F}$ as a subsheaf of $\mathscr{F}/a^n\mathscr{F}$ for any $n$, which means that the ``closures'' of $\mathscr{G}_1$ and $\mathscr{G}_2$ in $\mathscr{F}$ are the same.
More precisely, for any affine open $U=\Sph A$ of $X$, where $M=\Gamma(U,\mathscr{F})$, $N_i=\Gamma(U,\mathscr{G}_i)$ ($i=1,2$), the closures of $N_1$ and $N_2$ in $M$ coincide with each other.
Now by adhesiveness of $A$, $N_i$'s are closed in $M$, which implies $N_1=N_2$.
Hence $\mathscr{G}_1=\mathscr{G}_2$, as desired.

As a corollary, we have the following result.
\begin{cor}\label{cor-GHGA}
For proper and flat schemes over $\Spec V$, the map
$$
\Hom_{\Sch_V}(X,Y)\longrightarrow\Hom_{\Hs^{\ast}_V}(\het{X},\het{Y})
$$
$($where $\Hs^{\ast}_V$ denotes the category of henselian schemes adic over $\Sph V)$ is bijective.
\end{cor}

Indeed, by GFGA, the map is injective.
To show it is surjective, note that any morphism $\het{X}\rightarrow\het{Y}$ is represented by its graph $\Gamma\in\het{X}\times_V\het{Y}$, defined by a coherent ideal of $\het{X}\times_V\het{Y}=\het{(X\times_VY)}$.
Hence the claim follows from Theorem \ref{thm-subquotients}.

\section{Henselian rigid geometry}\label{sec-henselianrigidspaces}
\subsection{Henselian affinoid algebras}\label{sub-henselianaffinoidalgebra}
Let $V$ be an $a$-adically henselian valuation ring of height one, where $a\in\m_V\setminus\{0\}$, and $K=\Frac(V)$ the fractional field.
We write 
$$
{\textstyle K\{X_1,\ldots,X_n\}=V\{X_1,\ldots,X_n\}\otimes_VK\ (=V\{X_1,\ldots,X_n\}[\frac{1}{a}]),}
$$
and call it the {\em henselian Tate algebra} over $K$.
$K\{X_1,\ldots,X_n\}$ is a Noetherian $K$-algebra. 
\begin{dfn}\label{dfn-henselianaffinoidalg}{\rm 
A {\em henselian affinoid algebra} over $K$ is a $K$-algebra $\mathcal{A}$ of the form 
$$
\mathcal{A}=K\{X_1,\ldots,X_n\}/\mathfrak{a},
$$
where $\mathfrak{a}\subseteq K\{X_1,\ldots,X_n\}$ is an ideal.}
\end{dfn}

Since $K\{X_1,\ldots,X_n\}$ is Noetherian, there exists a finitely generated ideal $\til{\mathfrak{a}}\subseteq V\{X_1,\ldots,X_n\}$ such that $\mathfrak{a}=\til{\mathfrak{a}}[\frac{1}{a}]$.
One can replace $\til{\mathfrak{a}}$ by its $a$-saturation, which is still finitely generated. 
Hence there exists a $V$-flat model $A=V\{X_1,\ldots,X_n\}/\til{\mathfrak{a}}$ such that $A[\frac{1}{a}]\cong\mathcal{A}$.

It can be shown that henselian affinoid algebras are Jacobson (Nullstellensatz for henselian affinoid algebras), and have Noether normalization.

\subsection{Henselian rigid spaces}\label{sub-henselianrigidspaces}
In order to define rigid spaces from henselian schemes, we need the notion of admissible blow-ups, which can be defined quite similarly to the formal scheme case; it is defined locally as the henselization of the usual algebraic blow-ups along the admissible ideals, i.e., open finite type ideal sheaves.
Then the category $\CRh$ of coherent ($=$ quasi-compact and quasi-separated) henselian rigid spaces is the quotient category of the category $\CHs^{\ast}$ of coherent henselian schemes with adic morphisms, mod-out by all admissible blow-ups.
We have the quotient functor
$$
\rig\colon\CHs^{\ast}\longrightarrow\CRh,\quad X\longmapsto X^{\rig},
$$
which we call, similarly to the usual situation, the ``rig'' functor.

General henselian rigid spaces (not necessarily quasi-compact, nor quasi-separated) can be defined by means of ``birational patching'' similarly to the formal case as in \cite[Chap.\ {\bf II}, 2.2.(c)]{FK2}

One has the notion of {\em affinoids} in the similar way: A henselian rigid space $\mathscr{X}$ is an affinoid if it is isomorphic to a henselian rigid space of the form $(\Sph A)^{\rig}$.
If $A$ is a henselian of finite type algebra over $V$ (where $V$ is as in \S\ref{sub-henselianaffinoidalgebra}), then the affinoid $(\Sph A)^{\rig}$ is the finite type affinoid space over $K$ corresponding to the affinoid algebra $\mathcal{A}=A[\frac{1}{a}]$.

For a coherent henselian rigid space $\mathscr{X}=X^{\rig}$, the projective limit of all admissible blow-ups of $X$
$$
\ZR{\mathscr{X}}=\varprojlim_{X'\rightarrow X}X'
$$
taken in the category of topological spaces is a quasi-compact space (\cite[Chap.\ {\bf 0}, 2.2.(c)]{FK2}), and is independent of the choice of the henselian model $X$.
We call the topological space $\ZR{\mathscr{X}}$ the {\it Zariski-Riemann space} associated to $\mathscr{X}$ (cf.\ \cite[Chap.\ {\bf II}, 3.1]{FK2}).
It comes with the canonical sheaf of rings $\O^{\int}_{\mathscr{X}}$, the so-called {\em integral structure sheaf} (cf.\ \cite[Chap.\ {\bf II}, 3.2.(a)]{FK2}), such that the pair $(\ZR{\mathscr{X}},\O^{\int}_{\mathscr{X}})$ gives the projective limit of all admissible blow-ups $X'$ of $X$ in the category of locally ringed spaces.
The {\it rigid structure sheaf}, denoted by $\O_{\mathscr{X}}$, is constructed as
$$
\O_{\mathscr{X}}=\varinjlim_{n>0}\lHom_{\O^{\int}_{\mathscr{X}}}(\mathscr{I}_X^n\O^{\int}_{\mathscr{X}},\O^{\int}_{\mathscr{X}}),
$$
where $\mathscr{I}_X$ is an ideal of definition of $X$.
Similarly to the formal scheme situation, $(\ZR{\mathscr{X}},\O_{\mathscr{X}})$ is a locally ringed space (cf.\ \cite[Chap.\ {\bf II}, 3.2.(b)]{FK2}).
The Zariski-Riemann space $\ZR{\mathscr{X}}$, as well as the two structure sheaves, can be constructed by patching for general henselian rigid spaces.

Again, similarly to the formal schemes case, points of the Zariski-Riemann space $\ZR{\mathscr{X}}$ are classified by valuations (cf.\ \cite[Chap.\ {\bf II}, 3.3]{FK2}).
A {\em rigid point} on a henselian rigid space $\mathscr{X}$ is a morphism of henselian rigid spaces of the form $\alpha\colon(\Sph V)^{\rig}\rightarrow\mathscr{X}$, where $V$ is an $a$-adically henselian valuation ring (of arbitrary height), or equivalently, an adic map $\Sph V\rightarrow\ZR{\mathscr{X}}$ (where `adic' means the condition similar to that in \cite[Chap.\ {\bf II}, 3.3.(a)]{FK2}).
It follows that the points of $\ZR{\mathscr{X}}$ are in natural one to one correspondence with the equivalence class of rigid points, with the equivalence class generated by the relation by ``domination'' of valuation rings (cf.\ \cite[Chap.\ {\bf II}, 3.3.(a)]{FK2}).

\subsection{Classical points}\label{sub-classicalpoints}
The henselian rigid spaces under our consideration in the sequel are henselian rigid spaces of finite type over $K$.

A henselian rigid space $\mathscr{Z}$ is said to be {\em point-like} if it is coherent and reduced, having a unique minimal point in $\ZR{\mathscr{Z}}$.
\begin{dfn}\label{dfn-classicalpoints}{\rm 
Let $\mathscr{X}$ be a henselian rigid space of finite type over $K$.
A {\em classical point} of $\mathscr{X}$ is a point-like locally closed rigid subspace $\mathscr{Z}\subseteq\mathscr{X}$.}
\end{dfn}

It can be shown (cf.\ \cite[Chap.\ {\bf II}, (8.2.9)]{FK2}) that any classical point $\mathscr{Z}\hookrightarrow\mathscr{X}$ is a closed subspace.
Moreover, it is of the form $\mathscr{Z}=(\Sph W)^{\rig}$, where $W$ is finite, flat, and finitely presented over $V$, such that $s(\mathscr{Z})=\Spec W[\frac{1}{a}]$ is a point (cf.\ \cite[Chap.\ {\bf II}, (8.2.7)]{FK2}).
To show this and the following results, we use the following fact, which follows from the GHGA comparison for $\H^0$ (Theorem \ref{thm-H0}).
\begin{prop}\label{prop-GHGA2}
Let $A$ be a henselian finite type $V$-algebra, and $X\rightarrow\Sph A$ an admissible blow-up.
$($Note that $X\rightarrow\Sph A$ is the henselization of a blow-up of $\Spec A$.$)$
Then $A'=\Gamma(X,\O_X)$ is a finite $A$-algebra $($hence $a$-adically henselian and $a$-adically universally adhesive$)$ such that $A[\frac{1}{a}]\cong A'[\frac{1}{a}]$.
\end{prop}

\begin{prop}[{\rm cf.\ \cite[Chap.\ {\bf II}, (8.2.10)]{FK2}}]\label{prop-classicalpoint1}
Let $\mathscr{X}=(\Sph A)^{\rig}$ be a finite type affinoid over $K$, where $A$ is $V$-flat henselian finite type $V$-algebra.
For any classical point $\mathscr{Z}\hookrightarrow\mathscr{X}$, the image of the map $s(\mathscr{Z})\rightarrow s(\mathscr{X})=\Spec A[\frac{1}{a}]$ is a closed point.
This establishes a canonical bijection between the set of all classical points of $\mathscr{X}$ and the set of all closed points of the Noetherian scheme $s(\mathscr{X})$.
\hfill$\square$
\end{prop}

Let us remark here that, to show Proposition \ref{prop-classicalpoint1}, one needs Proposition \ref{prop-GHGA2}, for which the GHGA comparison for $\H^0$ (Theorem \ref{thm-H0}) is necessary.

We denote by $\mathscr{X}^{\cl}$ the set of all classical points of $\mathscr{X}$.
The formation $\mathscr{X}\mapsto\mathscr{X}^{\cl}$ is functorial (cf.\ \cite[Chap.\ {\bf II}, (8.2.14)]{FK2}).

\subsection{GAGA}\label{sub-GAGA}
Let $A$ be a henselian finite type $V$-algebra, and consider
$$
{\textstyle U=\Spec A[\frac{1}{a}]\hookrightarrow S=\Spec A}
$$
and the closed subset $D=\Spec A/aA$.

For a separated $U$-scheme $f\colon X\rightarrow U$ of finite type, we define the category $\Emb_{X|S}$ whose objects are the commutative diagrams 
$$
\xymatrix{X\ \ar@{^{(}->}[r]\ar[d]_f&\ovl{X}\ar[d]^{\ovl{f}}\\ U\ \ar@{^{(}->}[r]&S\rlap{,}}
$$
where $\ovl{f}\colon\ovl{X}\rightarrow S$ is a proper $S$-scheme, and $X\hookrightarrow\ovl{X}$ is a birational open immersion, with the $S$-morphisms $\ovl{X}\rightarrow\ovl{X}'$ that are $X$-admissible, i.e., isomorphisms on $X$.
The category $\Emb_{X|S}$ is cofiltered, and essentially small (cf.\ \cite[Chap.\ {\bf II}, (9.1.1)]{FK2}).
For any object given by the diagram above, we set
\begin{equation*}
\begin{split}
Z&=(\ovl{X}\times_SU)\setminus X,\\
\ovl{Z}&=\textrm{the closure of}\ Z\ \textrm{in}\ \ovl{X},\\
\til{X}&=\ovl{X}\setminus\ovl{Z}.
\end{split}
\end{equation*}
Consider the $a$-adic henselization $\het{\til{X}}\hookrightarrow\het{\ovl{X}}$ of the open immersion $\til{X}\hookrightarrow\ovl{X}$, which gives rise to an open immersion $(\het{\til{X}})^{\rig}\hookrightarrow(\het{\ovl{X}})^{\rig}$ of henselian rigid spaces.
Define 
$$
X^{\an}=\varinjlim(\het{\til{X}})^{\rig},
$$
where the inductive limit is taken along $\Emb^{\opp}_{X|S}$.
This gives rise to a functor, called the {\em GAGA functor}, from the category of separated of finite type $U$-schemes to the category of henselian rigid spaces locally of finite type over $K$.
It can be shown that GAGA functor commutes with fiber products (cf.\ \cite[Chap.\ {\bf II}, (9.1.10)]{FK2}).

As in the usual rigid geometry case, the GAGA theorems should follow from GHGA theorems, some of which we have mentioned in \S\ref{sub-GHGA}, which are sufficient for our later discussion.

\subsection{Quasi-finite morphisms}\label{sub-qfin}
Let $\varphi\colon \mathscr{X}\rightarrow\mathscr{Y}$ be a morphism between locally of finite type henselian rigid spaces over $K$.

\begin{dfn}\label{dfn-quasifinite}{\rm 
We say that the morphism $\varphi$ is {\em quasi-finite} if for any point $x\in\ZR{\mathscr{Y}}$, the fiber $\ZR{\varphi}^{-1}(x)$ is a finite set.}
\end{dfn}

Note that, if $x$ is a classical point, then $\ZR{\varphi}^{-1}(x)$ consists of finitely many classical points of $\mathscr{X}$.

A morphism $\varphi\colon\mathscr{X}\rightarrow\mathscr{Y}$ of coherent henselian rigid spaces is said to be {\em finite}, if it has a finite henselian model, i.e., there exists a finite morphism $f\colon X\rightarrow Y$ of henselian schemes such that $\varphi=f^{\rig}$.
(A morphism $f\colon X\rightarrow Y$ of henselian schemes is {\em finite} if for any affine open $V=\Sph B\subseteq Y$, $f^{-1}(V)$ is affine $\Sph A$ with $A$ finite over $B$.)
By the description of points of the associated Zariski-Riemann spaces in \S\ref{sub-henselianrigidspaces}, one sees the following: For a finite morphism $\varphi\colon\mathscr{X}\rightarrow\mathscr{Y}$ between locally of finite type henselian rigid spaces over $K$, and for any point $y\in\ZR{\mathscr{Y}}$, the fiber $\ZR{\varphi}^{-1}(y)$ is a finite set.
Hence, in particular, finite morphisms are quasi-finite.

\begin{rem}{\rm It will follow from the proof of our henselian version of ZMT (Theorem \ref{thm-ZMT} below) that a quasi-compact $\varphi\colon \mathscr{X}\rightarrow\mathscr{Y}$ is quasi-finite if and only if for any classical point $x$ of $\mathscr{Y}$ the fiber $\ZR{\varphi}^{-1}(x)$ is a finite set, since the proof only uses the finiteness of the fibers over classical points.
Indeed, to apply Theorem \ref{thm-ZMT}, one first reduce to the case where $\mathscr{X}$ and $\mathscr{Y}$ are affinoids, and if $\mathscr{Y}=(\Sph A)^{\rig}$, then take a finite type ring $B$ over $V$ such that $\het{B}\cong A$, and replace $\mathscr{Y}$ by $(\Spec B[\frac{1}{a}])^{\an}$.}
\end{rem}

\subsection{Zariski Main Theorem}\label{sub-ZMT}
Let $X$ be a separated of finite type scheme over $K$.
\begin{thm}\label{thm-ZMT}
Let $\varphi\colon\mathscr{U}\rightarrow X^{\an}$ be a quasi-finite $K$-map from a henselian affinoid space $\mathscr{U}$ of finite type over $K$.
Then there exists a finite morphism $g\colon W\rightarrow X$ with the diagram
$$
\xymatrix{\mathscr{U}\ar[r]_{j}\ar@/^1pc/[rr]^{\varphi}&W^{\an}\ar[r]_{g^{\an}}&X^{\an},}
$$
where $j$ is an open immersion.
\end{thm}

Let us sketch the proof.
Write $A$ as the inductive limit of the inductive system $\{A_{\lambda}\}$ of finite type $V$-algebras such that each transition map $A_{\lambda}\rightarrow A_{\mu}$ is \'etale with $A_{\lambda}/aA_{\lambda}\cong A_{\mu}/aA_{\mu}$, and that the henselization of each $A_{\lambda}$ coincides with $A$.
One can show, similarly to \cite[Chap.\ {\bf II}, \S9.2]{FK2}, that the given map $\varphi\colon\mathscr{U}=(\Sph A)^{\rig}\rightarrow X^{\an}$ comes from $\til{\varphi}\colon\Spec A[\frac{1}{a}]\rightarrow X$; this is where we have to use a GHGA theorem.

Here let us sketch the construction of $\til{\varphi}$.
In the notation as in \S\ref{sub-GAGA}, the given map $\varphi$ can extend to a map of henselian schemes $W\rightarrow\het{\til{X}}$ (for a suitable Nagata compactification $\ovl{X}$ of $X$), where $W$ is an admissible blow-up of $\Sph A$.
Notice that $W$ is the henselization of a blow-up $T\rightarrow\Spec A$ along an ideal supported in the closed fiber; i.e., $W$ is algebrizable $W=\het{T}$.
Now, by Corollary \ref{cor-GHGA}, one has a morphism $T\rightarrow\ovl{X}$ of schemes, which clearly gives $T\rightarrow\til{X}$.
Since $T\otimes_VV[\frac{1}{a}]\cong\Spec A[\frac{1}{a}]$, we get the desired $\til{\varphi}\colon\Spec A[\frac{1}{a}]\rightarrow X$ by passage to the generic fibers.

Now, there exists $\lambda$ such that $\til{\varphi}$ factors by $\til{\varphi}_{\lambda}\colon\Spec A_{\lambda}[\frac{1}{a}]\rightarrow X$, which one can show to be quasi-finite.
This follows from the classical Chevalley's Theorem \cite[$\mathbf{IV}$, (13.1.1)]{EGA} and the comparison of classical points and closed points as in Proposition \ref{prop-classicalpoint1}.
Then, since the original $\mathscr{U}=(\Sph A)^{\rig}$ is an affinoid subdomain of $(\Spec A_{\lambda}[\frac{1}{a}])^{\an}$, the desired theorem follows from the classical ZMT.

\begin{cor}\label{cor-ZMT}
Let $X$ be a separated finite type scheme over $K$, $\mathscr{U}$ a henselian rigid space of finite type over $K$, and $\varphi\colon \mathscr{U}\hookrightarrow X^{\an}$ an immersion.
Then there exists a {\em closed subscheme} $W\subseteq X$ that is smallest among those containing the image of $\mathscr{U}$ as an open subspace.
\end{cor}

Indeed, it suffices to show that there exists at least one closed $W$ containing $\mathscr{U}$ as an open subspace, and so we may assume $\mathscr{U}$ is an affinoid.
Then we apply ZMT to obtain $W$ finite over $X$, and replace it by the scheme-theoretic image (in the usual sense).

\subsection{Family of closed spaces}\label{sub-FCS}
\begin{thm}\label{thm-FCS}
Let $X$ be a proper scheme over $K$.
For any flat family $\mathscr{Y}\subseteq X^{\an}\times_K\mathscr{U}$ of closed subspaces in $X^{\an}$ over finite type henselian affinoid $\mathscr{U}=(\Sph A)^{\rig}$, there exists an affine finite type $K$-scheme $U$ such that
\begin{itemize}
\item[{\rm (a)}] $U^{\an}$ contains $\mathscr{U}$ as an affinoid subdomain$;$
\item[{\rm (b)}] $\mathscr{Y}$ extends to a flat family $Y\subseteq X\times_KU$ over $U$.
\end{itemize}
\end{thm}

To sketch the proof, write $A=\varinjlim_{\lambda}A_{\lambda}$ as before, and let $Y_{\lambda}$ be the scheme-theoretic image of $\mathscr{Y}$ in $X^{\an}\times_KU^{\an}_{\lambda}$, where $U_{\lambda}=\Spec A_{\lambda}[\frac{1}{a}]$.
One sees that, since $\mathscr{Y}\rightarrow\mathscr{U}$ is proper, the projective limit of $Y_{\lambda}\rightarrow U_{\lambda}$ recovers $\mathscr{Y}\rightarrow\mathscr{U}$ on passage to the analytification. 
By a standard limit arguement, there exists $\lambda$ such that $Y_{\lambda}\rightarrow U_{\lambda}$ is flat, which gives, therefore, a desired extension.

\frenchspacing
\begin{small}

\medskip\noindent
{\sc Department of Mathematics, Tokyo Institute of Technology, 2-12-1 Ookayama, Meguro, Tokyo 152-8551, Japan} (e-mail: {\tt bungen@math.titech.ac.jp})
\end{small}
\end{document}